\documentclass[11pt]{article}

\usepackage{bbm}
\usepackage[colorlinks=true,pagebackref,hyperindex]{hyperref}
\usepackage{mathrsfs} 
\usepackage[notcite,notref]{showkeys}
\usepackage[all]{xy}
\usepackage{color}
\usepackage{amsfonts}
\usepackage{amscd}

\definecolor{darkgreen}{rgb}{0.0, 0.7, 0.0}

\definecolor{cyan}{cmyk}{1,0,0,0}

\newcommand{\bdg}{\begin{dg}}
\newcommand{\edg}{\end{dg}}

\newtheorem{tm}{Theorem}[subsection]

\newtheorem{cor}[tm]{Corollary}

\newtheorem{??}[tm]{Question}

\newcommand{\ben}{\begin{enumerate}}
\newcommand{\een}{\end{enumerate}}
\newcommand{\bit}{\begin{itemize}}
\newcommand{\eit}{\end{itemize}}
\newcommand{\beq}{\begin{equation}}
\newcommand{\eeq}{\end{equation}}
\newcommand{\la}{\label}

\newcommand\ci{\cite}

\font\tenmsb=msbm10
\font\sevenmsb=msbm7
\font\fivemsb=msbm5
\newfam\msbfam
\textfont\msbfam=\tenmsb
\scriptfont\msbfam=\sevenmsb
\scriptscriptfont\msbfam=\fivemsb
\def\Bbb#1{{\fam\msbfam #1}}
\font\teneufm=eufm10
\font\seveneufm=eufm7
\font\fiveeufm=eufm5
\newfam\eufmfam
\textfont\eufmfam=\teneufm
\scriptfont\eufmfam=\seveneufm
\scriptscriptfont\eufmfam=\fiveeufm
\def\frak#1{{\fam\eufmfam\relax#1}}

\newcommand\rat{{\Bbb Q}}

\newcommand\comp{{\Bbb C}}

\newcommand{\bb}[1]{\Bbb{#1}}

\newcommand{\pc}[2]{ \,^{\frak p}\!{\cal H}^{#1}({#2})   }
\newcommand{\pcs}{ \,^{\frak p}\!{\cal H}   }

\title{Decomposition theorem for semi-simples}
\author{
Mark Andrea A.  de Cataldo\thanks{
Partially supported by N.S.F. grant DMS 1600515}}

\begin{document}

\maketitle

\begin{abstract}
We use standard constructions in algebraic geometry and homological algebra to 
extend the decomposition and hard Lefschetz theorems of
T. Mochizuki and C. Sabbah so that they remains valid without the quasi-projectivity
assumptions.
\end{abstract}

\tableofcontents

\section{Introduction}\la{intro}

 M. Kashiwara
\ci{Ka}  has put-forward a series of  conjectures
concerning  the behavior
of holonomic  semi-simple $D$-modules  on a complex algebraic variety
under proper push-forward and under taking
nearby/vanishing cycles. 

Inspired by this conjecture, T. Mochizuki \ci{Mo} has proved Kashiwara conjectures
in the very important case where one assumes  the holonomic $D$-modules
to be regular. Mochizuki's work built on earlier work by C. Sabbah \ci{Sa}.
Because of the regularity assumptions (see  \ci[p.2-3, Remark 6]{Sa}) for more context), part of their results can be expressed, via the Riemann-Hilbert correspondence, in the form of Theorem \ref{mosa} below. 

The methods employed in \ci{Mo,Sa} are  
essentially analytic. Moreover, \ci{Mo,Sa}  are placed  in the context of projective 
morphisms of quasi projective manifolds, so that Theorem \ref{mt} below, which generalizes  
Theorem \ref{mosa}, is not  directly affordable
by their methods: one would first need to extend aspects of their theory of polarizable pure  twistor $D$-modules  from projective manifolds
to complex algebraic varieties. To my knowledge, this extension is not in the literature.

V. Drinfeld \ci{Dr} has shown that an arithmetic conjecture by A. de Jong implies,
rather surprisingly and again under the regularity assumption,
 Kashiwara's conjectures. Drinfeld's proof uses also algebraic geometry  for varieties over finite fields.
 Note that \ci{Dr} allows 
for arbitrary characteristic-zero coefficients.
 de Jong's conjecture  has been proved by D. Gaitsgory \ci{Ga}  and by G. B\"ockle and C. Khare 
\ci{Bo-Ka}. 

The combination of the work in \ci{Dr,Ga,Bo-Ka} yields an arithmetic proof
of Theorems \ref{mosa} and of  \ref{mt} below.

The purpose of this note is to provide a proof of Theorem 
\ref{mt} that stems directly from Theorem \ref{mosa} and uses only simple reductions  based on standard constructions  in algebraic geometry.

{\bf Acknowledgments.} I am grateful to T. Mochizuki and to C. Sabbah for very  useful
remarks.

\section{Decomposition and relative hard Lefschetz for semi-simples}\la{s1}
\subsection{Statement}\la{stat}

A variety is a separated scheme of finite type over the field of complex numbers $\comp$.
For the necessary background  concerning what follows, the reader may consult
\ci{dCM}.
Given a variety $Y,$ we work with the rational and complex  constructible derived categories $D(Y,\rat)$ and  $D(Y,\comp)$
endowed with the middle-perversity $t$-structures, whose hearts, i.e. the respective categories of perverse sheaves
on $Y,$ are  denoted by $P(Y,\rat)$ and $P(Y,\comp)$, respectively.   The simple objects in $P(Y,\rat)$
and in $P(Y, \comp)$ have the form 
$IC_S(L),$ where $S$ is an  irreducible closed subvariety of $Y,$ $L$ is a simple
(i.e. irreducible)  complex/rational local system
defined on some dense open subset of the regular part of $S,$  and $IC$ stands for intersection complex.
 We say that $K\in D(Y,\rat)$
is semi-simple if it is isomorphic to the finite  direct sum of shifted simple perverse sheaves as above: 
$K\cong \oplus_b \pc{b}{K}[-b] \cong \oplus_b \oplus_{(S,L) \in EV_b} IC_S(L) [-b],$ where $\pcs^b$ denotes the $b$-th perverse cohomology sheaf functor, and  $EV_b$ is a uniquely determined finite set of pairs $(S,L)$ as above. Similarly, with $\comp$-coefficients.

Our starting point is the following result of T. Mochizuki \ci[\S14.5 and \S14.6]{Mo}, which generalizes one of C. Sabbah \ci{Sa}. In fact,
they both work in the  more refined setting of polarized pure twistor
$D$-modules and their results have immediate  and evident counterparts in the setting of the constructible derived category, which is the one of this note.

\begin{tm}\la{mosa}
Let $f:X \to Y$ be a projective map of irreducible quasi projective nonsingular  varieties. If $K\in P(X,\comp)$ is semi-simple,
then $f_* K \in D(Y,\comp)$ is semi-simple. The relative hard Lefschetz theorem holds.
\end{tm}

Even if the methods in \ci{Mo} seem to require the smoothness and quasi projectivity 
assumptions, as well as $\comp$-coefficients, one can deduce   the following more general statement.
We have nothing to say concerning  the refined context of polarizable pure twistor $D$-modules.

\begin{tm}\la{mt}
Let $f: X \to Y$ be a proper map of varieties. If $K\in P(X,\rat)$ is semi-simple, then $f_*K \in D(Y,\rat)$
is semi-simple. If $f$ is projective, then the relative hard Lefschetz theorem holds.
\end{tm}

We first show how to deduce  the $D(Y,\comp)$-version of Theorem \ref{mt} from Theorem \ref{mosa}. Then 
we show how the $D(Y,\comp)$-version implies formally the $D(Y,\rat)$-version. 

The reader should have no difficulty in replacing $\rat$ with any field of characteristic zero
and proving the same result.

\subsection{Proof of Theorem \ref{mt} for $D(Y,\comp)$}

Theorem \ref{mosa} is stated for  $\comp$-coefficients. In this section, we use this statement to deduce
Theorem \ref{mt} for $\comp$-coefficients, i.e. to deduce Corollary \ref{mtc} below.

The theorem will be reduced to several special cases, where  we progressively relax the hypotheses
on $f$, from 
projective,  to quasi projective,  to proper,  and on $X$ and $Y$, from smooth quasi projective, to
quasi projective, to arbitrary. These conditions will be denoted symbolically by $(f_{proj}, X^{sm}_{qp}, \ldots)$.
For example,
we summarize the hypotheses of Theorem \ref{mosa} graphically as follows: 
\[(f_{proj}, X^{sm}_{qp}, Y^{sm}_{qp}) \qquad \mbox{
($f$ projective,  $X$ and $Y$ smooth and quasi projective)}.
\]

Our goal is to establish Corollary \ref{mtc}  as an immediate consequence of the five  following claims.

\ben
\item
{\em Theorem \ref{mosa} holds for $(f_{proj}, X^{sm}_{qp}, Y_{qp})$.} 

Choose any closed 
embedding $g:Y \to 
\bb{U}$  of $Y$ into a Zariski-dense open subvariety
$\bb{U} \subseteq \bb{P}$ of some projective space.
Apply Theorem \ref{mosa} to $h:=g\circ f$
and observe that, modulo the natural identification of  the objects in $D(Y,\comp)$
with the  ones in $D(\bb{U},\comp)$ supported on $Y,$  we have $h_*K=f_*K$.  

\item
{\em Theorem \ref{mosa} holds for $(f_{proj}, X_{qp}, Y_{qp})$.}

 Pick a resolution of the singularities  $g:Z \to X$ of $X$
with $g$ projective. Let $X^o \subseteq X_{reg}
\subseteq X$ be a dense Zariski open subset on which  the simple local system $M$ is defined and over which $g$ is an isomorphism.
Let $IC_Z(M) \in P(Z,\comp)$ be the  intersection complex  on $Z$ with coefficients  in the local system $M$ transplanted to $g^{-1}(U^o)$.
Apply 1. to $g$ and $h$. Observe that $IC_X(M)$ is a direct summand
of $g_* IC_Z(M)$. Deduce that  $f_* IC_X(M)$ is a direct summand of $h_* IC_Z(M)$ so that the first part of Theorem \ref{mosa}
holds for $(f_{proj}, X_{qp}, Y_{qp})$. In order to prove the second part of Theorem \ref{mosa}, i.e.  the relative hard Lefschetz theorem for $f$, we argue as in 
\ci{dCM}, Lemma 5.1.1: we do not need self-duality to conclude: the argument gives injectivity;
by dualizing we get surjectivity for the dual of the hard Lefschetz maps; this dualized map is the hard Lefschetz map 
for $f,$ $IC_X(M)^\vee$ and the  $f$-ample $\eta \in H^2(X,\comp)$;
by switching the roles of $M$ and $M^\vee,$ we see that the relative hard Lefschetz theorem maps are isomorphisms.
(N.B.: we may impose self-duality artificially, by replacing $M$ with $M\oplus M^{\vee}$ and reach the same conclusion.)

\item
{\em Theorem \ref{mosa} holds for $(f_{proj}, X_{qp}, Y)$.}

 Let $Y=\cup_i Y_i$ be an affine open covering. Let $f_i: X_i:= f^{-1}(Y_i) \to Y_i$
be the obvious maps. By 2., the relative Hard Lefschetz holds for $f_i.$ Since the relative hard Lefschetz
maps are defined over $Y$ and they are isomorphisms over the $Y_i,$ the relative hard Lefschetz holds for $f$ over $Y$.
By the Deligne-Lefschetz criterion \ci{De}, we have $f_*K \cong  \oplus_b \pc{b}{f_*K}[-b]$. It remains to show
that the  $P^b:=\pc{b}{f_*K}$ are semi-simple. By 2., the $P^b_{|Y_i}$  are semi-simple after restriction to the open affine  $Y_i.$
By  a repeated use of  the the splitting criterion \ci{dCM}, Lemma 4.1.3\footnote{let $P$ be a perverse sheaf on a variety $Z$; let $Z=U\coprod Z$ be Whitney-stratified  in such a way that
$U\subseteq Z$ is open  and union of strata,  $S\subseteq Z$ is a  closed stratum, and $P$ is  cohomologically constructible with respect to the stratification; Lemma 4.1.3 in  \ci{dCM} is an iff criterion for the splitting of $P$ into the intermediate extension $j_{!*} (P_{|U})$  to $Z$ of the restriction $P_{|U}$ of $P$ to $U$,  direct sum a local system on $S$ placed in cohomological degree minus the codimension of the stratum; the criterion is local in the classical and even in the Zariski topology}
applied in the context of a
Whitney stratification of $Y$ w.r.t. which the $P^b$ are cohomologically constructible, we deduce that the $P^b$ split as direct sum of intersection 
complexes with coefficients in some local systems. (Note that  \ci{dCM}, Assumption 4.1.1  is fulfilled in view of \ci{dCM}, Remark 4.1.2, because we already know that $P_b$ splits as desired over the open $Y_i$.)
We need to verify that these  local systems are semi-simple. Since a local system
on an integral  normal variety is semisimple if and only if it is semisimple after restriction
to a Zariski dense open subvariety, the desired semi-simplicity  can be checked by restriction to the chosen affine covering 
of $Y,$ where we can apply 2.

\item
{\em Theorem \ref{mosa} holds for $(f_{proj}, X, Y).$}

 As it was pointed out in 3., the relative hard Lefschetz can be verified  on an affine covering  $Y=\cup_i Y_i.$ The resulting $X_i$ are then quasi-projective and we can apply 3. For the semisimplicity of the direct image 
$f_* IC_X(M),$ we take a Chow envelope $g:Z \to X$  of $X$ ($Z$ quasi projective, $g$ projective and birational); we produce $IC_Z(M)$ as above and we deduce the semisimplicity of $f_* IC_X(M)$ from the one --established in 3.-- of $h_*IC_Z(M),$ as it was done in 2.

\item 
{\em The semisimplicity statement in Theorem \ref{mosa} holds for $(f_{proper}, X, Y)$.}

 Take a Chow envelope
$g: Z \to X$ of $f$ ($g$ birational, $g$ and $h:= f\circ g$ projective).  Produce $IC_Z(M)$ as above.
Apply   4.  and deduce that $f_* IC_X(M)$
is a direct summand of the semi-simple $h_* IC_Z (M).$

\een

The above, together with the obvious remark that it is enough to prove Theorem \ref{mt} 
in the case when $X,Y$ are irreducible and $K=IC_X(M),$ yields the following 

\begin{cor}\la{mtc} Theorem \ref{mt} holds for $\comp$-coefficients.
\end{cor}

\subsection{Theorem \ref{mt} for $D(Y,\comp)$ implies the same for $D(Y,\rat)$}\la{rt5}

Let $f$ be projective. Then we have the relative hard Lefschetz for $\comp$-coefficients, hence for $\rat$-coefficients as well.
By the Deligne-Lefschetz criterion, we have the isomorphism $f_* K \cong \oplus_b \pc{b}{f_*K} [-b]$ in $D(Y,\rat)$.
We need to show that each $P^b:=\pc{b}{f_*K} [-b]$ is semi-simple in $P(Y,\rat)$.
Note that extending the coefficients from $\rat$ to $\comp$ is a $t$-exact
functor $D(Y, \rat) \to D(Y, \comp)$. In particular, the formation of $P^b$ is compatible with  complexification.
By arguing as in point 3. of the previous section, we see
that each $P^b$ is a direct sum of intersection complexes $IC_S(L)$, where  the $L$ are rational local systems (note that [dCM], Assumption 4.1.1 is now  fulfilled in view of [dCM], Remark 4.1.2, because we already know that the complexification of $P^b$ splits as desired over $Y$).
We need to verify that  each $L$ is a semi-simple rational 
local system. We know its complexification is, hence so is $L$, in fact:
let $0\to L'\to L\to L'' \to 0$ be an extension of rational locally constant
sheaves on $S^o$;
it is classified by an element $e \in H^1(S^o, {L''}^*\otimes L')$; this element becomes  trivial after complexification, hence it is trivial over $\rat$.

If $f$ is proper, we take a Chow envelope $g: Z \to X$ of $f,$ we set $h:=f\circ g$ and we deduce semisimplicity of $f_*$ from the semisimplicity
of $h_*$  ($h$ is projective) as in point 5. of the previous section.

\providecommand{\bysame}{\leavevmode \hbox \o3em
{\hrulefill}\thinspace}

Mark Andrea A. de Cataldo, Department of Mathematics,
Stony Brook University,
Stony Brook,  NY 11794, USA;
mark.decataldo@stonybrook.edu

\end{document}